%

\magnification=1200
\baselineskip18pt

\def\ms{\medskip}

\def\phi{\varphi}

\def\su{\subseteq}
\def\a{\alpha}
\def\b{\beta}

\def\e{\nu}
\def\d{\delta}
\def\l{\lambda}
\def\k{\kappa}
\def\z{\zeta}

\def\lng{\langle}
\def\rng{\rangle}
\def\ov{\overline}
\def\sm{\setminus}


\def\acc{{\rm acc}\,}

\def\cf{{\rm cf}}

\def\otp{{\rm otp}\,}
\def\ran{{\rm  ran}}


\def\proof{\smallbreak\noindent{\sl Proof}: }

\newcount\itemno
\def\itm{\advance\itemno1 \item{(\number\itemno)}}
\def\ritm{\advance\itemno1 \item{)\number\itemno(}}
\def\startitm{\itemno=-1 }
\def\startaitm{\itemno=0}
\def\aitm{\advance\itemno1
\item{(\letter\itemno)}}

\def\letter#1{\ifcase#1 \or a\or b\or c\or d\or e\or f\or g\or h\or
i\or j\or k\or l\or m\or n\or o\or p\or q\or r\or s\or t\or u\or v\or
w\or x\or y\or z\else\toomanyconditions\fi}

\def\imply{\Rightarrow}


\def\rest{\mathop |}
\centerline{\bf A PROOF OF SHELAH'S PARTITION THEOREM}

\centerline{Menachem Kojman, Carnegie-Mellon University}

\bigskip

The following is  self contained presentation of Shelah's recent  
proof
 of the partition relation  $\Big({\buildrel {\mu^+} \over \mu}\Big)
 \longrightarrow
\Big({\buildrel {\mu+1}\over \mu}\Big)^{1,1}_{<\cf\mu}$
 for a  singular strong limit $\mu$ violating the GCH. The notation
 $\Big({\buildrel {\mu^+} \over \mu}\Big) \longrightarrow
\Big({\buildrel {\mu+1}\over \mu}\Big)^{1,1}_{<\cf\mu}$
 means: for every coloring $c$ of $\mu^+\times \mu$ by less than
 $\cf\mu$ many colors there are $A\su\mu^+$ with $\otp A=\mu+1$ and
 $B\su \mu$ with $\otp B=\mu$ such that $c$ is constant on $A\times  
B$.

 The proof here is re-arranged slightly differently than the proof in
 the forthcoming [Sh 513] so that no use of other results of Shelah  
is
 made, except for Fact 2 below, which comes from pcf theory.  In  
other
 words, we avoid here using the ideal $I[\l]$ from [Sh 420] and the
 tools from [Sh 108]; now it is not that reading those two papers is  
a
 bad idea --- on the contrary, I have been intending to do so myself  
for
 a number of years now. It is only that the proof is accessible
 directly.

 The pcf theory needed to obtain 2 below will be available also in a
 survey paper [K] on pcf and $I[\l]$.

\ms\noindent
{\bf 1. Theorem}: Suppose that $\mu$ is a strong limit singular  
cardinal
 and $2^\mu>\mu^+$. Then

$$\Big({\buildrel {\mu^+} \over \mu}\Big) \longrightarrow
\Big({\buildrel {\mu+1}\over \mu}\Big)^{1,1}_{<\cf\mu}$$

\proof We prove actually a stronger claim: for every function
 $c:(\mu^+\times \mu)\to \theta$, for $\theta<\cf \mu$, there are  
$A\su
 \mu^+$ and $B\su \mu$ with $\otp A=\mu+1$ and $\otp B=\mu$ such that
 the fuction $c\rest (A\times B)$ does not depend on the first
 coordinate. This clearly implies the theorem.

Let $\k$ denote $\cf \mu$. Fix an increasing sequence $\ov \mu=\lng
 \mu_i: i<\k\rng$ cofinal in $\mu$ such that $\mu_0>\k$.
The assumptions we made about $\mu$ imply the following:

\ms\noindent
{\bf 2.  Pcf Fact}: There is an increasing sequence of regular  
cardinals
 $\ov \l=\lng \l_i: i<\k\rng$ with $\sup \ov \l=\mu$ such that
 $\Pi_i\l_i/{J^{bd}_\k}$ is $\mu^{++}$-directed, where $J^{bd}_\k$ is
 the ideal of bounded subsets of $\k$.

This Fact follows from $\mu$ being a strong limit and $2^\mu>\mu^+$  
via
 pcf theory. For details see chapter 8 of Shelah's book or [K].

We may thin out $\ov \l$ and assume that $\l_i>2^{\mu_i^+}$.

Suppose that $c:(\mu^+\times \mu)\to \theta$ is given for some
 $\theta<\k$. We need to produce $A$ and $B$ as above. These sets  
will
 be constructed in $\k$ many approximations, after some preparation.

\ms

Fix a function $F$ from $[\mu^+]^2$ to $\k$ that such that for all
 $i<\k$ the set $a^\a_i:=\{\b<\a:F(\a,\b)\le i\}$ has cardinality at
 most $\mu_i^+$.  Thus $\a=\bigcup_{i<\k}a^\a_i$.

Let $\chi$ be a sufficiently large regular cardinal. We define by  
double
 induction of $\mu^+\times \k$ a matrix $\{M_{\a,i}:\a<\mu^+,i<\k\}$  
of
 elementary submodels of $(H(\chi),\in)$, satisfying:

\startitm

\itm $M_{\a,i}\prec (H(\chi,\in)$, $||M_{\a,i}||=2^{\mu_i^+}$ and
 ${}^{\mu_i^+}M_{\a,i}\su M_{\a,i}$ ($M_{\a,i}$ is closed under
 sequences of length $\mu_i^+$).

\itm $\a,c,\ov \mu,\ov \l$ and $F$ belong to $M_{\a,i}$ and
 $\{M_{\b,j}:(\b,j)<_{lx}(\a,i)\}$ belongs to $M_{\a,i}$.

\ms

There is no problem to choose $M_{\a,i}$ so that is satisfies the
 conditions above.

We make a few simple observations about this array or models:

\ms\noindent
{\bf 3. Fact}:
\startitm
\itm If $(\b,j)<_{lx}(\a,i)$ and $\b\in M_{\a,i}$ then $M_{\b,j}\in
 M_{\a,i}$.

\itm If $M_{\b,j}\in M_{\a,i}$ and $j\le i$ then $M_{\b,j}\su  
M_{\a,i}$
 and hence $M_{\b,j}\prec M_{\a,i}$.

\itm $M_{\a,j}\prec M_{\a,i}$ for all $\a<\mu^+$ and $j<i<\k$.

\itm $\a\su \bigcup_iM_{\a,i}$ for all $\a<\mu^+$

\itm For all $\b<\a<\mu^+$ for an end segment of $i<\k$ it holds that
 $M_{\b,i}\su M_{\a,i}$ and hence
$M_{\b,i}\prec M_{\a,i}$.

\proof
Clause (0) is follows from the demand that
 $\{M_{\b,j}:(\b,j)<_{lx}(\a,i)\}\in M_{\a,i}$ and the fact that  
$\k\su
 \mu_i\su M_{\a,i}$, so $i\in M_{\a,i}$, and therefore $M_{\b,j}$ is
 definable from parameters in $M_{\a,i}$. Being an elementary  
submodel,
 $M_{\a,i}$ contains every set definable from parameters in  
$M_{\a,i}$.

To see clause (1) suppose that $M_{\b,j}\in M_{\a,i}$ and that $j\le  
i$.
 By elementarity of $M_{\a,i}$ there is a bijection
 $\varphi:2^{\mu_i^+}\to M_{\b,j}$ in $M_{\a,i}$. As $2^{\mu_i^+}\su
 M_{\a,i}$, also $\ran \varphi\su M_{\a,i}$ and hence $M_{\b,j}\su
 M_{\a,i}$. Since also $M_{\b,j}\prec (H(\chi),\in)$ and  
$M_{\a,i}\prec
 (H(\chi),\in)$, necessarily $M_{\b,j}\prec M_{\a,i}$ and (1) holds.

Clause (2) follows from the previous two and the fact that $\a\in
 M_{\a,i}$.

To prove (3) use the fact that $a^\a_i\in M_{\a,i}$ and also  
$a^\a_i\su
 M_{\a,i}$ for all $i<\k$. Therefore for all $i\ge F(\a,b)$ it holds
 that $\b\in M_{\a,i}$. Thus (3) holds.

The last clause follows from the previous ones.

A conclusion of those facts is the following:

\ms\noindent
{\bf 4.  Fact}: The sequence $\ov M_\a=\lng M_{\a,i}:i<\k\rng$ is
 increasing in $\prec$, $\a\su \bigcup _iM_{\a,i}$ and if  
$\b<\a<\mu^+$
 then $\ov M_\b\in_{J^{bd}_\k}\ov M_\a$, $\ov M_b\su_{J^{bd}_\k}\ov
 M_\a$, and even $\ov M_b\prec_{J^{bd}_\k}\ov M_\a$, namely for all
 sufficiently large $i<\k$ we have that $M_{\b,i}\in M_{\a,i}$,
 $M_{\b,i}\su M_{\a,i}$ and $M_{\b,i}\prec M_{\a,i}$.

\ms

For every $\a<\mu^+$  and $i<\k$ define $f_\a(i)=\sup M_{\a,i}\cap
 \lambda_i$. As we assumed that $\l_i>2^{\mu_i^+}=||M_{\a,i}||$, it
 follows by the regularity of $\l_i$ that $f_\a(i)\in \l_i$, for all
 $i<\k$ and therefore $f_\a\in \Pi \l_i$ for all $\a<\mu^+$.

Furthermore, if $\b<\a<\mu^+$ then from some $i_{\a,\b}<\k$ onwards
 $M_{\b,i}\in M_{\a,i}$ and therefore (as $\ov \l\su  M_{\a,i}$)
 $f_\b(i)\in M_{\a,i}$ and hence $f_\b(i)<f_\a(i)$ on an end segment  
of
 $\k$, or $f_\b<_{J^{bd}_\k}f_\a$. Thus $\ov f=\lng  
f_\a:\a<\mu^+\rng$
 is increasing in $< _{J^{bd}_\k}$.

Use Fact 2 above to find a bound $f^*\in \Pi \l_i$ to $\ov f$ in $\le
 _{J^{bd}_\k}$.

Using $f^*$ and the coloring $c$, define $g_\a(i)=c(\a,f^*(i))$ for  
all
 $\a<\mu^+$ and $i<\k$.
The function $g_\a$ specifies the $c$-type of $\a$ over the sequence
 $\lng f^*(i):i<\k\rng$.

As there are only $\theta^\k<\mu^+=\cf\mu^+$ many possible such  
types,
 we find a function $g^*:\k\to \theta$ so that $A:=\{\a<\mu^+:
 g_\a=g^*\}$ is unbounded in $\mu^+$.

\ms

Let us find now
by induction on $\zeta<\mu^+$ an increasing continuous chain of
 elementary submodels $\ov N=\lng N_\a:\zeta<\mu^+\rng$ satisfying:

\startitm

\itm $\mu\su N_\zeta\prec (H(\chi,\epsilon)$ and $||N_\zeta||=\mu$

\itm  $A$, $g^*$ and $\{M_{\a,i}:\a<\mu^+,i<\k\}$ belong to $N_0$

\ms
Let $E=\{\zeta<\mu^+:\zeta=N_\zeta\cap \mu^+\}$. This is a club of
 $\mu^+$.

By induction on $i<\k$ we choose a strictly increasing sequence of
 ordinals $\d_i<\mu^+$ satisfying:

\startaitm
\aitm $\d_i\in \acc E$ (that is, $\d_i$ is an accumulation point of  
$E$)
 and

\aitm $\cf \d_i=\mu_i^+$.

\ms
Observe that $\d_i>\sup \{\d_\e:\e<i\}$ for all $i<\k$, because
 $\cf\d_i=\mu_i^+$. This enables us to choose $\a(i)\in \d_i\sm
 \sup\{\d_\e : \e<i\}$ for every  $i<\k$.

We also observe that if $\a\in N_\zeta$ then $M_{\a,i}\prec N_\z$ for
 $i<\k$. Therefore, if $\zeta\in E$, then $M_{\a,i}\prec N_\zeta$ for
 all $\a<\z$ and $i<\k$.

Pick $\a(*)\in A\sm \sup\{\d_i:i<\k\}$.

We define now by induction on $i<\k$ sets $A_i,B_i$ and an index
 $j(i)<\k$ such that the following conditions hold:

\startaitm

\aitm $j(i)>i$ and $i_1<i_2\imply \l_{j(i_1)}<\mu_{j(i_2)}$

\aitm For any two ordinals $\sigma<\tau$ in the set $\{\d_\e:\e\le
 i\}\cup \{\a_\e:\e\le i\}\cup \{\a(*)\}$ it holds that $\ov M_\sigma
 \prec \ov M_\tau$ and $f_\sigma < f_\tau$ on the end segment
 $(j(i),\k)$ of $\k$.

\aitm $A_i\su A\cap \d_i$,  $\otp A_i= \mu_i^+$ and $A_i\in
 M_{\d_i,j(i)}$.

\aitm $B_i\su \l_{j(i)}\sm \sup\{\l_{j(\e)}: \e<i\}$, $\otp
 B_i=\l_{j(i)}$ and $B_i\in M_{\d_i,j(B_i)}$ for some $j(B_i)<\k$.  
Also,
 $B_\e\in M_{\d_i,j(i)}$ for all $\e<i$.

\aitm If $\a\in \bigcup_{\e\le i}A_i\cup\{\a(*)\}$ and $\b\in B_\e$  
for
 some $\e\le i$ then $c(\a,\b)=g^*(j(\e))$.

\ms

If the induction is carried out successfully, then by (e)   it  
follows
 that if $\a\in A=\bigcup_{i<\k}A_i\cup\{\a(*)\}$ and $\b\in
 B=\bigcup_{i<\k}B_i$ then $c(\a,\b)=g^*(j(i))$ for the (unique)  
first
 $i$ satisfying $\l_{j(i)}>\b$.  From  (c) and (d) it follows that  
$\otp
 A=\mu+1$ and $\otp B=\mu$. Thus $A,B$ are as required by the  
theorem.

Suppose, then, that $A_\e,B_\e$ and $j(\e)$ are defined for all  
$\e<i$
 and satisfy the conditions above.

Since $\a(i)>\e$ for every $\e<i$, there is some $j(\e)<\k$ such that
 $B_\e,A_\e,j(\e)\in M_{\a(i),j}$ for $j\ge j(\e)$. Let $j_0<\k$ be
 large enough so that $B_\e,A_\e,j(\e)\in M_{\a(i),j_0}$ for all $\e<
i$ and so that $\mu_{j_0} > \l_{j(\e)}$ for all $\e<i$. This can be  
done
 as there are less than $\k$ many $\e$-s.

We have, then, $B_\e\in M_{\a(i),j_0}$ for all $\e<i$ or
 $\{B_\e:\e<i\}\su M_{\a(i),j_0}$. As $M_{\a(i),j_0}$ is closed under
 sequences of length at most $\mu_{j_0}^+>\k$ we also have that $\lng
 B_\e:\e<i\rng\in M_{\a(i),j_0}$.
Similarly, $\lng A_\e:\e<i\rng\in M_{\a(i),j_0}$ and $\lng j(\e)
 :\e<i\rng\in M_{\a(i),j_0}$.

Since $\d_i$ is an accumulation point of $E$ and has cofinality
 $\mu_i^+$, we can find an increasing  sequence $\lng \zeta_ \epsilon  
:
 \epsilon<\mu_i^+\rng$ of elements of $E$ with $\zeta_0>\a(i)$.

For every $\zeta_\epsilon$ in the sequence we chose, $\a(i)\in
 \zeta_\epsilon\su N_{\zeta_\epsilon}$, and therefore
 $M_{\a(i),j_0}\prec N_{\zeta_\epsilon}$ and hence $\lng B_\e:  
\e<i\rng,
 \lng j(\e):\e<i\rng\in N_{\zeta_\epsilon}$.

 For every $\epsilon <\mu_i^+$ the ordinal $\a(*)$ satisfies  in
 $(H(\chi_,\in)$  the following formula $\varphi(x,\zeta_\epsilon)$
 (when substituted for $x$):

$$\varphi(x,\zeta_\epsilon):= x\in A\;\&\;x>\zeta_\epsilon\;\&\;
 (\forall \e<i)(\b\in B_\e\imply c(x,\b)=g^*(j(\e)))
\leqno{(1)}$$

Since all the parameters in this sentence --- namely $A$,  $\lng  
B_\e:
 \e<i\rng$, $\lng j(\e):\e<i\rng$, $c$, $g^*$ and  $\zeta_\epsilon$  
---
 belong to $N_{\zeta_{\epsilon+1}}$ and the latter is an elementary
 submodel of $(H(\chi),\in)$, there is an ordinal $\gamma_\epsilon\in
 N_{\zeta_{\epsilon+1}}$ such that
 $\varphi(\gamma_\epsilon,\zeta_\epsilon)$ holds. Clearly,
 $\zeta_\epsilon<\gamma_\epsilon<\zeta_{\epsilon+1}<\d_i$.

Let $A'_i:=\{\gamma_{\epsilon+1}: \epsilon <\mu_i^+\}$. We have shown
 that $A'_i\su A\cap (\a(i),\d_i)$ and every $\a\in A'_i$ satisfies  
that
 $c(\a,\b)=g^*(j(i))$ for the first $i$ such that $\l_{j(i)}>\b$.   
Each
 member of $A'_i$  belongs to $M_{\d_i,j}$ for some $j<\k$, since
 $\d_i\su \bigcup _{j<\k}M_{\d_i,j}$. Because $\mu_i^+>\k$ is  
regular,
 there must be some index $j_1<\k$ such that $A(i)=A'(i)\cap
 M_{\d_i,j_1}$ has cardinality $\mu_i^+$. Let $A(i)$ be the set   
$A_i$
 we need to define. This takes care of the first two parts in (c).

Let $j(i) \ge \max\{j_1,j_0\}$ be large enough so  that $A_i\in
 M_{\d_i,j_1}$ and $M_{\d_i,j(i)}\prec M_{\a(*),j(i)}$, and also such
 that $f_{\d_i}(j(i))<f^*(j(i))$. Now the remaining part of (c), (a)  
and
 (b) are also satisfied.

Work now in $M_{\a(*),j(i)}$. We know that $\lng A_\e:\e<i\rng,
 A_i,\a(*)\in M_{\a(*),j(i)}$ and that also the function $\e\mapsto
 j(\e)$ for $\e<\i$ belongs to $M_{\a(*),j(i)}$, because all  
functions
 from $\k$ to $\k$ belong to it.

 Therefore the following set is definable in $M_{\a(*),j(i)}$:

$$B:=\{\b<\l_{j(i)}: c(\a,\b)=g^*(j(i)) {\rm \;for\;all\;} \a\in
 \bigcup\limits_{\e\le i}A_\e\cup\{\a(*)\}\}\leqno{(2)}$$

Observe that $f^*(j(i))$ belongs to the set $B$ defined in (2)  
because
 $\bigcup\limits_{\e\le i}A_\e\cup\{\a(*)\}\su A$, but that since
 $f^*(j(i))>f_{\d_i}(j(i))=\sup M_{\d_i,j(i)}\cap \l_{j(i)}$ it does  
not
 belong to $M_{\d_i,j(i))}$. This shows that $B$ has no bound in
 $M_{\d_i,j(i)}\cap \l_{j(i)}$. We conclude, then, that $B$ is  
unbounded
 below $\l_{j(i)}$: being definable in $M_{\d_i,j(i)}$, if there were  
a
 bound to $B$ below $\l_{j(i)}$ there would be one in  
$M_{\d_i,j(i)}$;
 but there is not.

Using the same argument as before, we find some $j(B)<\k$ such that
 $B_i=B\cap M_{\d_i,j(B)}\sm \sup \{ \l_{j(\e)} : \e <i\}$ belongs to
 $M_{\d_i, j(B)}$ and has cardinality $\l_{j(i)}$. Now (d) and (e)  
are
 also satisfied.

This completes the induction, and the proof as well.

\bigbreak

\noindent

\noindent
{\bf References}

\bigskip
\noindent
[K] M.~ Kojman, {\sl PCF and $I[\l]$}, in preparation

\ms
\noindent
[S] S.~Shelah, {\bf Cardinal Arithmetic}, Oxford University Press,  
1994
 Oxford.

\ms
\noindent
[Sh 108] S.~ Shelah, {\sl On successors of singular cardinal} in  
Logic
 Colloquium 78, volume 97 of Stud. Logic. Foundations Math, pp.
 357--380, North Holland 1979

\ms
\noindent
[Sh 420] S. Shelah, {\sl Advances in Cardinal Arithmetic} in  
Proceedings
 of the Banff Conference in algebra.

\ms
\noindent
[Sh 513] S.~ Shelah {\sl pcf and infinite free subsets in an  
algebra},
 in preparation.

\end